\documentclass{IEEEtran}

\usepackage{amsmath,amsfonts,cite,times,url,hyperref,bm,graphicx}

\newtheorem{example}{Example}

\def\BibTeX{{\rm B\kern-.05em{\sc i\kern-.025em b}\kern-.08em
    T\kern-.1667em\lower.7ex\hbox{E}\kern-.125emX}}
\begin{document}

\title{Convex Pollution Control of Wastewater Treatment Systems}

\author{Joshua A. Taylor, \IEEEmembership{Senior Member, IEEE}
\thanks{Josh A. Taylor is with the Department of Electrical and Computer Engineering,
       New Jersey Institute of Technology, Newark, NJ, USA. E-mail:
        {\tt\small jat94@njit.edu}}
}

\maketitle

\begin{abstract}

We design a model-predictive controller for managing the actuators in sewer networks. It minimizes flooding and combined-sewer overflow during rain and pollution at other times. To make the problem tractable, we use a convex relaxation of the microbial growth kinetics and a physically motivated linearization of the mass flow bilinearities. With these approximations, the trajectory optimization in each control period is a second-order cone program. In simulation, the controller releases roughly 15\% less pollutant mass than a conventional controller while treating nearly the same volume of flow. It does so by better balancing the flow over the treatment plants and over time.

\end{abstract}

\begin{IEEEkeywords}
Wastewater treatment, pollution control, model-predictive control, second-order cone programming, convex relaxation
\end{IEEEkeywords}

\section{Introduction}\label{sec:intro}

Sewer networks can retain wastewater in storage tanks and long pipes. They can control the flow of wastewater with actuators like gates, valves, and pumps. Together, this storage capacity and actuation enable system operators to manage the flow arriving at the treatment plants.

Today, this flexibility is most often used during rain to mitigate flooding and combined-sewer overflow (CSO)~\cite{pleau2005globalfull}. At other times, the system is well within capacity, and the flexibility can be leveraged to improve treatment efficiency, thereby reducing pollution. In this paper, we design a model-predictive controller (MPC) that handles both tasks---it minimizes flooding and CSO during high flow, and minimizes pollutant release at other times.

Controllers that solely mitigate flooding and CSO are known as volume-based~\cite{ocampo2010model,ocampo2013applicationfull,garcia2015modeling}. They do not account for pollutant concentrations or treatment processes, and as such are based on tractable, linear models. Immission-based controllers account for the sewer network, the treatment plants, and the receiving environment~\cite{schutze1999optimisation}. While they have high theoretical performance, they are difficult to implement. The middle ground is known as pollution-based control, which seeks to maximize the pollutant mass arriving at the plants, but without modeling the treatment plants or receiving environment~\cite{mahmoodian2017pollution,sun2020integrated}. In all cases, it is common to use MPC, in which the controller solves a trajectory optimization in each time period~\cite{rawlings2017model}.

Our controller is most nearly pollution-based, differing only in that we also account for the biochemical treatment processes in the plants. The model in our controller's trajectory optimization has two nonconvexities: microbial growth kinetics (e.g., with the Monod~\cite{monod1949growth} or Contois~\cite{contois1959kinetics} rates), and mass flows with bilinear products of flow rates and concentrations.

In Section~\ref{sec:socapprox}, we convexify the growth kinetics using the second-order cone (SOC) relaxations of~\cite{taylor2021grad,taylor2022cob}. This renders the trajectory optimization biconvex. In prior work~\cite{taylor2024predictive}, we used the alternating direction method of multipliers (ADMM)~\cite{boyd2011distributed} to solve this biconvex optimization as a sequence of convex problems. This performs well, but has two drawbacks: ADMM is not guaranteed to converge on biconvex problems, and solving a sequence of ten to twenty convex optimizations is time-consuming.

In Section~\ref{sec:concapprox}, we eliminate the bilinearities with the following novel approximations.
\begin{itemize}
\item We only model concentrations and reaction rates in the plants because plant volumes are roughly constant.
\item We only model flow rates and volumes in the sewer network because concentrations do not change significantly in the pipes. We estimate plant inlet concentrations from influent forecasts and forward simulation of the system.
\end{itemize}
After these approximations, the trajectory optimization is a second-order cone program (SOCP), which can be reliably solved at large scales and well within real-time~\cite{Boyd1998SOCP,boyd2004convex}. This is the first pollution-based MPC that relies entirely on convex optimization.

In Section~\ref{sec:simulations}, we present a simulation-based case study in which we compare our new controller with conventional volume-based MPC. With no algorithmic specialization or tuning, both controllers compute decisions between one and two orders of magnitude faster than real-time. Our controller releases roughly 15\% less pollutant mass while treating nearly the same volume of sewage as the volume-based controller. It achieves this reduction by better balancing the flow over the treatment plants and over time, leading to higher overall reaction rates and lower effluent concentrations in the plants.

\section{Modeling}\label{sec:modeling}

\subsection{Sewer network}\label{sec:sewer}
The system consists of $s$ tanks interconnected by flows. Denote the set of tanks $\mathcal{S}$. There are three types of tanks:
\begin{itemize}
\item treatment plants, $\mathcal{T}$,
\item real tanks, $\mathcal{R}$, and
\item virtual tanks, $\mathcal{V}$, which are long pipes with the ability to store wastewater.
\end{itemize}

$V_i$ is the volume of wastewater in tank $i\in\mathcal{S}$. The flow rates into and out of tank $i$ from and to outside the system are $Q_i^{\textrm{in}}$ and $Q_i^{\textrm{out}}$. They could respectively represent, e.g., influent from streets and buildings, and the effluent from a treatment plant. Let $\mathcal{P}\subset\mathcal{S}\times\mathcal{S}$ be the set of pipes. Let $Q_{ij}$ denote the flow rate from tank $i$ to tank $j$ if $ij\in\mathcal{P}$.

Some of the pipe flow rates are controlled by actuators such as pumps, detention gates, and diversion gates. We denote the set of actuators $\mathcal{A}\subset \mathcal{P}$. We model the actuators in Section~\ref{sec:act}.

Let $Q^\textrm{F}_i$ denote flood flow rate out of tank $i$. Only virtual tanks can flood, so that $Q^\textrm{F}_i=0$ if $i\notin\mathcal{V}$. Denote the total flow rates to and from tank $i$ by
\begin{subequations}
\begin{align}
\tilde{Q}_i^{\textrm{in}} &= Q_i^{\textrm{in}}+\sum_{k\in{\mathcal{S}}}Q_{ki}\\
\tilde{Q}_i^{\textrm{out}} &= Q^\textrm{F}_i + Q_i^{\textrm{out}}+\sum_{k\in{\mathcal{S}}}Q_{ik}.
\end{align}

The volume in non-plant tanks evolves as
\label{dynamics}
\begin{align}
\frac{dV_i}{dt}=\tilde{Q}_i^{\textrm{in}} - \tilde{Q}_i^{\textrm{out}} ,\quad i\in\mathcal{S}\setminus\mathcal{T}.\label{dVdt}
\end{align}
The non-plant volumes satisfy the constraint
\[
0\leq V_i\leq V_i^{\max},\quad i\in\mathcal{S}\setminus\mathcal{T}.
\]
For $i\in\mathcal{V}$, if $V_i$ would exceed $V_i^{\max}$, $Q^\textrm{F}_i$ is the flow rate that makes it equal to $V_i^{\max}$.

A plant only receives inflow from other tanks and only has one outflow. The volume of water in the treatment plants is constant, i.e., $V_i=\bar{V}_i$. CSO occurs if there is flow above a maximum rate, i.e., if
\begin{align}
Q^\textrm{CSO}_i=\max\left\{\tilde{Q}_i^{\textrm{in}}-Q_i^{\max},0\right\}
\end{align}
is greater than zero for some $i\in\mathcal{T}$. This implies that
\begin{align}
\tilde{Q}_i^{\textrm{in}}- Q^\textrm{CSO}_i&= Q_i^{\textrm{out}},\quad i\in\mathcal{T}.
\end{align}

\subsection{Biochemical processes}\label{sec:bio}

We follow the setup in Section 1.5 of~\cite{bastin2013line}. There are $m$ substrates and biomasses in each perfectly mixed tank. $\xi_i\in\mathbb{R}^m_+$ is the process state vector of tank $i\in\mathcal{S}$, which contains the concentrations of the substrates and biomasses. $\xi_i^{\textrm{in}}\in\mathbb{R}^m_+$ is the influent concentration vector of tank $i\in\mathcal{S}$.

Biochemical reactions occur in the treatment plants. There are $r_i$ different types of reactions that convert substrates to other substrates and biomasses in plant $i\in\mathcal{T}$. $\phi_i(\xi_i)\in\mathbb{R}^{r_i}_+$ is a vector of the reaction kinetics in plant $i$. The elements of $\phi_i(\xi_i)$ are concave functions that can be relaxed to SOC constraints, as we discuss in Section~\ref{sec:socapprox}. We show how to do this for Monod and Contois kinetics later in Examples \ref{ex:contois} and \ref{ex:monod}.

Let $\kappa_i\in\mathbb{R}^{m\times r}$ be the stoichiometric matrix relating the reaction vector, $\phi_i(\xi_i)$, to the evolution of the process state in plant $i$. Note that while all tanks contain the same types of concentrations, they can each have different reactions and stoichiometry. We refer the reader to~\cite{taylor2022cob} and~\cite{bastin2013line} for examples.

The dynamics of the concentrations are derived by taking a mass balance in each tank. The dynamics in non-plant tank $i\in\mathcal{S}\setminus\mathcal{T}$ are
\begin{align}
\frac{dV_{i}\xi_i}{dt} &=V_i\frac{d\xi_i}{dt} + \frac{dV_i}{dt}\xi_i\nonumber\\
&= Q_{i}^{\textrm{in}}\xi_i^{\textrm{in}}-\tilde{Q}_{i}^{\textrm{out}}\xi_i+\sum_{j\in\mathcal{S}}Q_{ji}\xi_j.\label{dynscalarV}
\end{align}

The inlet concentration of plant $i\in\mathcal{T}$ is
\begin{align}
\xi^{\textrm{in}}_i&=\frac{\sum_{j\in\mathcal{S}}Q_{ji}\xi_j}{\tilde{Q}_{i}^{\textrm{in}}}.\label{eq:xiin}
\end{align}
Because the plants have constant volumes, we can write their dynamics
\begin{align}
\frac{d\xi_i}{dt} &=\kappa_i\phi_i(\xi_i) +\frac{Q_{i}^{\textrm{out}}}{\bar{V}_{i}}\left(\xi^{\textrm{in}}_i-\xi_i\right),\quad i\in\mathcal{T}.\label{dynscalarP}
\end{align}
\end{subequations}

The complete dynamics are given by (\ref{dynamics}), with initial conditions for the volumes $V$, and concentrations, $\xi$. They are nonlinear due to the mass flow terms and microbial growth kinetics in the latter equations.

\subsection{Discretization and time delays}\label{sec:disc}
To make (\ref{dynamics}) compatible with finite-dimensional optimization, we must approximate it with a discrete time system. We use the Adams-Moulton implicit linear multistep scheme~\cite{butcher2016numerical} for its higher order accuracy and regular step sizes, which are convenient when sensor measurements and influent forecast data are at regular time intervals. Note that the first-order Adams-Moulton scheme is the trapezoidal method.

Index the time periods $n\in\mathcal{N}=\{1,...,\tau\}$, and let $\delta$ be the length of each time period. $V(n)$ and $\xi(n)$ are the vectors of volumes and concentrations at time $n\delta$. For $n\in\mathcal{N}$, define
\begin{subequations}
\begin{align}
&\mathcal{D}_{V_i}(n) = \tilde{Q}_i^{\textrm{in}}(n) - \tilde{Q}_i^{\textrm{out}}(n),\; i\in\mathcal{S}\setminus\mathcal{T}\\
&\mathcal{D}_{\xi_i}(n) = \kappa_i\phi_i(\xi_i(n)) +\frac{Q_{i}^{\textrm{out}}(n)}{\bar{V}_{i}}\left(\xi^{\textrm{in}}_i(n)-\xi_i(n)\right),\; i\in\mathcal{T}\label{eq:Dxi}\\
&\mathcal{D}_{M_i}(n) = Q_{i}^{\textrm{in}}(n)\xi_i^{\textrm{in}}(n) -\tilde{Q}_{i}^{\textrm{out}}(n)\xi_i(n) \nonumber\\
&\hspace{20mm}+\sum_{j\in\mathcal{S}}Q_{ji}(n)\xi_j(n),\; i\in\mathcal{S}\setminus\mathcal{T}.
\end{align}
\end{subequations}
The $K^{\textrm{th}}$ order discretization of (\ref{dynamics}) at time $n\in\mathcal{N}$ is 
\begin{subequations}
\label{eq:Moulton}
\begin{align}
&V_i(n)-V_i(n-1) = \delta\sum_{k=0}^K \alpha_k\mathcal{D}_{V_i}(n-k),\; i\in\mathcal{S}\setminus\mathcal{T}\label{discdynV}\\
&\xi_i(n)-\xi_i(n-1) = \delta\sum_{k=0}^K \alpha_k\mathcal{D}_{\xi_i}(n-k),\; i\in\mathcal{T}\label{discdynxi}\\
&V_i(n)\xi_i(n)-V_i(n-1)\xi_i(n-1) =\nonumber\\
&\hspace{20mm} \delta\sum_{k=0}^K \alpha_k\mathcal{D}_{M_i}(n-k),\; i\in\mathcal{S}\setminus\mathcal{T}.\label{discdynM}
\end{align}
\end{subequations}
The $\alpha_k$ are coefficients that can be found in, e.g.,~\cite{butcher2016numerical}. Equations (\ref{discdynV}) and (\ref{discdynxi}) will be constraints in the controller's trajectory optimization in Section~\ref{sec:optimization}.

Time delays occur due to transit in the pipes. As a result, the concentration arriving at one tank is the concentration departing another tank at an earlier time. This can be accounted for in (\ref{eq:Moulton}) by subtracting the delays in the argument of $\mathcal{D}_{V_i}(\cdot)$, $\mathcal{D}_{\xi_i}(\cdot)$, and $\mathcal{D}_{M_i}(\cdot)$.

Let $\tau^{\textrm{IC}}$ be the maximum over the largest time delay, the order of the Adams-Moulton scheme, and three. The initial conditions are the values of the state up to $\tau^{\textrm{IC}}$ periods before zero, i.e., $V^{\textrm{IC}}(l)$ and $\xi^{\textrm{IC}}(l)$, $l=-\tau^{\textrm{IC}},...,-1$. Here $\tau^{\textrm{IC}}$ has been chosen large enough to capture time delays, the earlier periods used by the discretization scheme, and earlier values of actuated flows later used in the slope and curvature objectives.

\subsection{Pipes and actuation}\label{sec:act}

The flows through the pipes depend on the actuators and the volumes in the tanks. This section describes the constraints on the volumes and flow rates due to the pipe network and actuators. We denote the set described by the constraints in this section by $\Omega$.

\begin{itemize}

\item The flow rates through the treatment plants satisfy
\[
Q_i^{\min }\leq Q_i^{\textrm{out}} \leq Q_i^{\max},\quad i\in\mathcal{T}.
\]

\item For all pipes $ij\in\mathcal{P}$,
\[
0\leq Q_{ij}(n)\leq Q_{ij}^{\max},
\]
and for all volumes $i\in\mathcal{S}\setminus\mathcal{P}$,
\[
0\leq V_{i}(n)\leq V_{i}^{\max}.
\]

\item Some pipes, $\mathcal{U}\in\mathcal{P}$, receive uncontrolled flow from upstream tanks. For $ij\in\mathcal{U}$, we approximate the flow rate as a linear function of the tank's volume:
\[
Q_{ij}(n) = \beta_iV_i(n),
\]
where $\beta_i$ is a constant~\cite{ocampo2010model}.

\item The flow rate through a pump or detention gate, $ij\in\mathcal{A}$, is either limited to a range,
\[
Q_{ij}^{\min}\leq Q_{ij}(n)\leq Q_{ij}^{\max},
\]
or limited by the stored volume in an upstream tank,
\[
Q_{ij}(n) \leq \beta_iV_i(n).
\]
\item Diversion gates are modeled as real tanks with zero volume, i.e.,
\[
V_i(n)=0,\quad i\in\mathcal{J}\subset\mathcal{R}.
\]
This implies that the flow rates in and out sum to zero at all times, and that we can choose the flow rates in the outgoing pipes subject to this constraint. The outgoing pipes are in the set of actuators, $\mathcal{A}$.
\item The flow through an actuator $ij\in\mathcal{A}$ might be held constant over every $\Delta$ minutes, e.g., if the time step is $\delta=5$ minutes and actuator setpoints are updated every $\Delta=$ 15 minutes. For each $n$ such that $n \textrm{ mod } \Delta/\delta=0$, we have
\[
Q_{ij}(n) = Q_{ij}(n+m), \; m=1,...,\Delta/\delta-1.
\]
\end{itemize}

\subsection{Performance metrics}\label{sec:metrics}

All of the metrics are linear or quadratic except for pollutant release, which is bilinear. A linear approximation for use in the optimization objective is given in Section~\ref{sec:concapprox}.

\begin{itemize}
\item \textit{Flooding.} Total volume of wastewater that exits the system before the treatment plants:
\[
\delta\sum_{n\in\mathcal{N}}\sum_{i\in\mathcal{V}}Q^{\textrm{F}}_i(n).
\]
\item \textit{Combined sewer overflow.} Total volume of wastewater that bypasses the treatment plants:
\[
\sum_{n\in\mathcal{N}}\sum_{i\in\mathcal{T}}Q^{\textrm{CSO}}_i(n).
\]
\item \textit{Pollutant release.} Total pollutant mass released from the plants into the receiving environment:
\[
w_{\textrm{PR}}^{\top}\sum_{n\in\mathcal{N}}\sum_{i\in\mathcal{T}}Q_i^{\textrm{out}}(n)\xi_i(n),
\]
where $w_{\textrm{PR}}$ is a vector of weights assigning relative priority to each substance.
\item \textit{Regulation violation.} Plant effluent concentrations in excess of regulatory limits:
\[
w_{\textrm{RV}}^{\top}\sum_{n\in\mathcal{N}}\sum_{i\in\mathcal{T}}\max\left\{
\xi_i(n) - \xi_i^{\max},0
\right\}.
\]
\item \textit{Smoothness of actuation.} Overly frequent or rapid changes to the setpoint of an actuator can cause wear and tear. The change in actuation can be quantified by slope and curvature, which are written
\[
\sum_{n\in\mathcal{N}}\sum_{ij\in\mathcal{A}} (Q_{ij}(n)-Q_{ij}(n-1))^2,
\]
and
\[
\sum_{n\in\mathcal{N}}\sum_{ij\in\mathcal{A}} (Q_{ij}(n)-2Q_{ij}(n-1)+Q_{ij}(n-2))^2.
\]
\end{itemize}

The following additional metrics, when used in the optimization objective, can improve controller performance. 
\begin{itemize}
\item \textit{Microbial growth.} Total microbial growth:
\[
w_{\textrm{MG}}^{\top}\sum_{n\in\mathcal{N}}\sum_{i\in\mathcal{T}}\bar{V}_i(n) T(n).
\]
\item \textit{Final and total volume.} Let $\tau$ be the last time period. The final volume in the system is
\[
\sum_{i\in\mathcal{S}\notin\mathcal{T}}V_i(\tau).
\]
This can prevent a receding horizon controller from restricting flow to lower pollutant release when the time horizon is short. The total volume in the system is
\[
\sum_{n\in\mathcal{N}}\sum_{i\in\mathcal{S}}V_i(n).
\]
This promotes increased flow through the system so as to keep the tanks empty.
\item \textit{Balanced plant utilization and outflow.} This pushes all plants to use the same fraction of their capacity.
\[
\sum_{n\in\mathcal{N}}\sum_{i\in\mathcal{T}} \left(\frac{Q_i^{\textrm{out}}(n)}{Q_i^{\textrm{max}}} - \frac{\sum_{i\in\mathcal{T}}Q_i^{\textrm{out}}(n)}{\sum_{i\in\mathcal{T}}Q_i^{\textrm{max}}}\right)^2.
\]
\item \textit{Balanced outflow.} This evens the flow through the plants over time.
\[
\sum_{n\in\mathcal{N}}\sum_{i\in\mathcal{T}}\left(\frac{\tau Q_i^{\textrm{out}}(n) - \sum_{n\in\mathcal{N}}Q_i^{\textrm{out}}(n)}{Q_i^{\textrm{max}}}\right)^2.
\]
\end{itemize}

\section{Controller}

We design an MPC~\cite{rawlings2017model} based on the model in Section~\ref{sec:modeling}. Consider a sequence of time periods, each of duration $\Delta$. Let $\tau^{\textrm{IC}}$ be as in the end of Section~\ref{sec:disc}.

We denote observations (or estimates) of the system state variables with a breve, e.g., $\breve{V}_i(n)$. The state observation in period $n$ consists of
\begin{itemize}
\item the volumes in all non-plant tanks, $\breve{V}_i(n-l)$, $i\in\mathcal{S}\setminus\mathcal{T}$,
\item the concentrations in the treatment plants, $\breve{\xi}_i(n-l)$, $i\in\mathcal{T}$, amd
\item the actuator setpoints, $\breve{Q}_{ij}(n-l)$, $ij\in\mathcal{A}$,
\end{itemize}
each for $l=1,...\tau^{\textrm{IC}}$.

The control in period $n$ is $Q_{ij}(n)$, $ij\in\mathcal{A}$, and $Q_{i}(n)$, $ij\in\mathcal{T}$ the flow rate setpoints of all actuators. The actuators then take appropriate action to realize the flow rate, e.g., pump speed or detention gate position. The control is computed by solving the optimization problem $\texttt{Traj}(n)$, which depends on information from up to period $n-1$. $\texttt{Traj}(\cdot)$ is an SOCP that we describe in Section~\ref{sec:optimization}.

The following sequence of events occurs in a given time period $n$.
\begin{enumerate}
\item At the start of period $n$, implement the control, $Q_{ij}(n)$, $ij\in\mathcal{A}$ and $Q_{i}(n)$, $ij\in\mathcal{T}$. This was computed during the previous period, $n-1$, by solving  $\texttt{Traj}(n)$.
\item The controller observes the state in period $n$. 
\item If desired, update the parameters of $\texttt{Traj}(\cdot)$, e.g., influent forecasts.
\item Compute control for the next period, $Q_{ij}(n+1)$, $ij\in\mathcal{A}$, by solving  \texttt{Traj}$(n+1)$.
\item Move forward in time to the next period, $n+1$, and go to Step 1.
\end{enumerate}

\subsection{Trajectory optimization}\label{sec:optimization}

This section describes $\texttt{Traj}(\cdot)$, the optimization problem the controller solves in each period. The full discretized dynamics, (\ref{eq:Moulton}), are difficult to optimize over because of the nonlinear microbial growth kinetics and the bilinear mass flows, both of which lead to nonconvex constraints. The pollutant release term in the objective is also bilinear. The intractability of such a problem would make real-time implementation haphazard.

We obtain a tractable SOCP by making the following two approximations.
\begin{itemize}
\item SOC relaxation of the microbial growth kinetics (Section~\ref{sec:socapprox}).
\item Only modeling concentrations and reaction rates in the plants, and only modeling volumes and flow rates in the network (Section~\ref{sec:concapprox}). This also requires approximation of the plants' inlet and outlet concentrations.
\end{itemize}
Section~\ref{sec:trajectory} states $\texttt{Traj}(\cdot)$. Let $\mathcal{N}_{\textrm{T}}=\{0,...,H\}$ be the sequence of time periods $\texttt{Traj}(\cdot)$ optimizes over.

\subsubsection{Convex relaxation of microbial growth kinetics}\label{sec:socapprox}

The term $\phi_i(\xi_i(n))$ in (\ref{eq:Dxi}) is nonlinear and makes the corresponding constraint, (\ref{discdynxi}), nonconvex. We relax it as follows. Introduce vectors of variables $T_i(n)$ and set it equal to the growth kinetics:
\begin{align}
T_i(n)=\phi_i(\xi_i(n)),\quad i\in\mathcal{T}, n\in\mathcal{N}_{\textrm{T}}.\label{eq:Tphi}
\end{align}
By adding this constraint, we can replace the growth rates in (\ref{discdynxi}) with $T_i(n)$. Now replace (\ref{eq:Tphi}) with the inequality
\begin{align}
T_i(n)\leq\phi_i(\xi_i(n)),\quad i\in\mathcal{T}, n\in\mathcal{N}_{\textrm{T}}.\label{eq:TphiR}
\end{align}
Because we only consider concave growth kinetics, (\ref{eq:TphiR}) is a convex constraint. It is a relaxation because all values that are feasible for (\ref{eq:Tphi}) are feasible for (\ref{eq:TphiR}), but not vice versa. This relaxation was first introduced and analyzed in~\cite{taylor2021grad}. It is often exact, meaning the optimal solution satisfies (\ref{eq:TphiR}) with equality. If not, exactness can be limited by adding a linear underestimator, as described in~\cite{taylor2021grad}.

The inequality (\ref{eq:TphiR}) can be written in SOC form for the growth rates in the following cases; other examples are given in~\cite{taylor2022cob}. For clarity, we omit the period index, $n$. In each case, there is a single scalar substrate concentration, $S$, biomass concentration, $X$, and growth kinetics, $T$.

\begin{example}[Contois]\label{ex:contois}
If the growth rate is Contois~\cite{contois1959kinetics}, constraint (\ref{eq:Tphi}) takes the form
\[
T\leq \frac{\mu SX}{k_{\textrm{C}}X+S},
\]
where $\mu$ and $k_{\textrm{C}}$ are the maximum growth rate and half-saturation constant. As shown in~\cite{taylor2021grad}, this can be written as the SOC constraint
\begin{align}
\left\|\left[\begin{array}{c}
\mu S\\
k_{\textrm{C}}T\\
\mu k_{\textrm{C}}X 
\end{array}\right]\right\|\leq \mu k_{\textrm{C}}X +  \mu S - k_{\textrm{C}}T.\label{eq:SOCC}
\end{align}
\end{example}

\begin{example}[Monod with constant biomass]\label{ex:monod}
If the growth rate is Monod~\cite{monod1949growth}, constraint (\ref{eq:Tphi}) takes the form
\[
T\leq \frac{\mu SX}{k_{\textrm{M}}+S}.
\]
The right-hand side is not concave. It becomes concave if we assume that the biomass in each time period is not an optimization variable, but an exogenous parameter, i.e., $X=\bar{X}$. This approximation is often valid because the biomass concentration in the treatment plants is orders of magnitude larger and slower-varying than the substrate concentrations. In this case, as shown in~\cite{taylor2021grad}, it can be written as the SOC constraint
\begin{align}
\left\|\left[\begin{array}{c}
\mu S \bar{X}\\
k_{\textrm{M}}T\\
\mu k_{\textrm{M}}\bar{X} 
\end{array}\right]\right\|\leq \mu k_{\textrm{M}}\bar{X} +  \mu S\bar{X} - k_{\textrm{M}}T.\label{eq:SOCM}
\end{align}
\end{example}

\begin{figure}[h]
\centering
\includegraphics[width=\columnwidth]{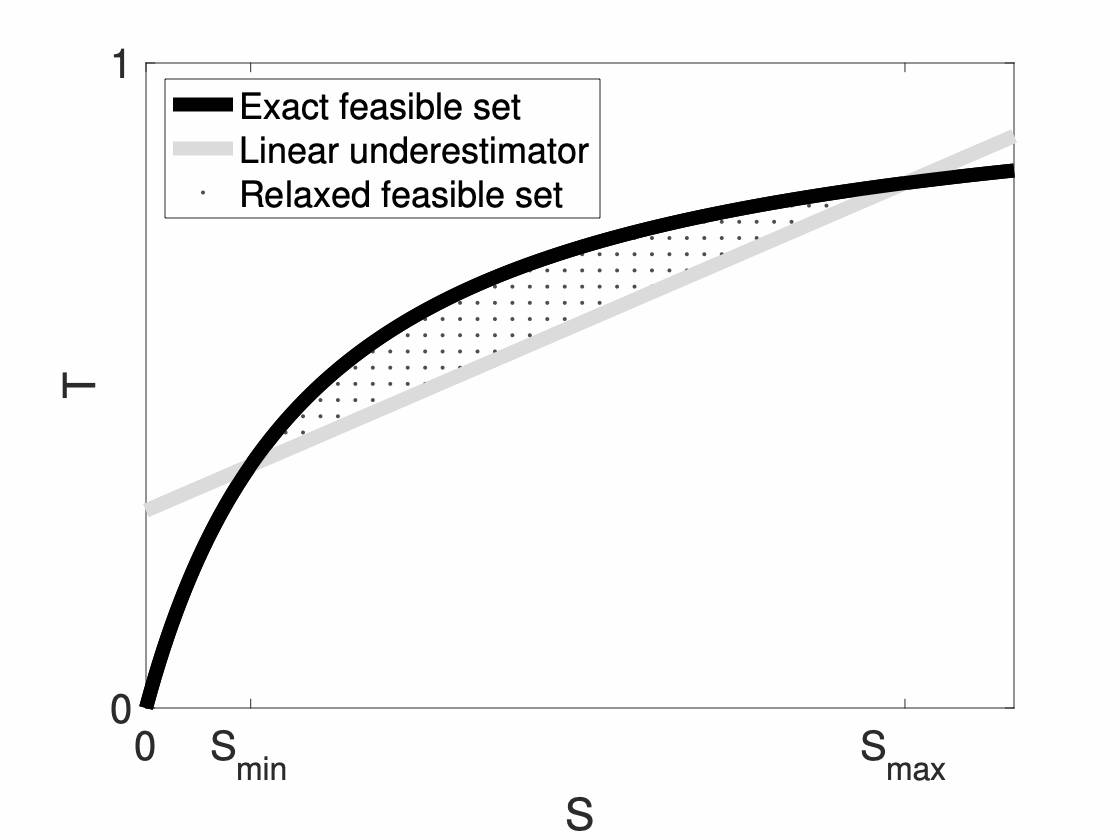}
\caption{The feasible sets described by the Monod growth rate and its second-order cone relaxation, with a linear underestimator.}
\label{fig:relaxationMonodunder}
\end{figure}

Figure~\ref{fig:relaxationMonodunder} shows the feasible set described by Constraint (\ref{eq:SOCM}), with a linear underestimator to limit potential inexactness. $S_{\min}$ and $S_{\max}$ are minimum and maximum substrate values, which can be zero and a large number if not known a priori.

%

\subsubsection{Linearization of bilinearities}\label{sec:concapprox}

We linearize the bilinearities in the discretized dynamics, (\ref{eq:Moulton}), as follows. First, omit (\ref{discdynM}). This eliminates the mass flow bilinearities in the pipe network. It is a reasonable approximation because the concentrations in the network do not evolve significantly or appear in the performance metrics.

Equation (\ref{discdynxi}) has two bilinearities because $\xi^{\textrm{in}}_i(n)$, as defined in (\ref{eq:xiin}), and $\xi_i(n)$ multiply the variable, $Q_i^{\textrm{out}}(n)$. To remove the bilinearities, we replace the variables $\xi^{\textrm{in}}_i(n)$ and $\xi_i(n)$ with the parameters $\hat{\xi}^{\textrm{in}}_i(n)$ and $\hat{\xi}_i(n)$, which are estimates of the inlet concentration and concentration in plant $i\in\mathcal{T}$ in period $n$. These can be obtained in several ways, including simulation of (\ref{dynamics}) with nominal future values of the controls taken from the solution of the trajectory optimization in the previous control period.

Define
\begin{align*}
\hat{\mathcal{D}}_{\xi_i}(n) &= \kappa_iT_i(n) +\frac{Q_{i}^{\textrm{out}}(n)}{\bar{V}_{i}}\left(\hat{\xi}^{\textrm{in}}_i(n)-\hat{\xi}_i(n)\right)
\end{align*}
for $i\in\mathcal{T}, n\in\mathcal{N}_{\textrm{T}}$. We approximate (\ref{discdynxi}) as
\begin{align}
\xi_i(n)-\xi_i(n-1) &= \delta\sum_{k=0}^K \alpha_k\hat{\mathcal{D}}_{\xi_i}(n-k)\label{discdynxi1}.
\end{align}

The pollutant release metric in Section~\ref{sec:metrics} is bilinear. In the objective of $\texttt{Traj}(\cdot)$, we replace this term with the linear approximation
\begin{align}
w_{\textrm{PR}}^{\top}\sum_{n\in\mathcal{N}_{\textrm{T}}}\sum_{i\in\mathcal{T}}Q_i^{\textrm{out}}(n)\hat{\xi}_i(n).\label{eq:obj:prlin}
\end{align}

\subsubsection{Second-order cone program}\label{sec:trajectory}

The objective, which we denote $\mathcal{F}(Q,V,\xi,T)$, is a weighted sum of the performance metrics in Section~\ref{sec:metrics}, summed over the time periods in $\mathcal{N}_{\textrm{T}}$. With the linear approximation of pollutant release, (\ref{eq:obj:prlin}), it contains only linear and convex quadratic terms. The SOCP $\texttt{Traj}(\cdot)$ is given below.

\begin{subequations}
\label{trajectory}
\begin{align}
\min_{Q,V,\xi,T}\quad&\mathcal{F}(Q,V,\xi,T) \label{P0}\\
\textrm{s.t.}\quad&\left(Q,V\right)\in\Omega\label{P7}\\
&\textrm{for }i\in\mathcal{T},l\in\mathcal{N}_{\textrm{T}}:\nonumber\\
&\hspace{5mm}\xi_i(l)\geq0, T_i(l)\geq0\label{P0}\\
&\hspace{5mm} T_i(l)\leq\phi_i(\xi_i(l))\label{P1}\\
&\textrm{for }i\in\mathcal{T}, l\in\mathcal{N}_{\textrm{T}}:\nonumber\\
&\hspace{5mm} \xi_i(l)-\xi_i(l-1) = \delta\sum_{k=0}^K \alpha_k\hat{\mathcal{D}}_{\xi_i}(l-k) \label{P2}\\
&\textrm{for }i\in\mathcal{S}\setminus\mathcal{T}, l\in\mathcal{N}_{\textrm{T}}:\nonumber\\
&\hspace{5mm} V_i(l)-V_i(l-1) = \delta\sum_{k=0}^K \alpha_k\mathcal{D}_{V_i}(l-k)\label{P3}\\
&\textrm{for }l=-\tau^{\textrm{IC}},...,-1:\nonumber\\
&\hspace{5mm}Q_i(l) = Q^{\textrm{IC}}_{i}(l),\quad i\in\mathcal{A} \label{P4}\\
&\hspace{5mm}V_i(l) = V^{\textrm{IC}}_{i}(l),\quad i\in\mathcal{S}\setminus\mathcal{T} \label{P5}\\
&\hspace{5mm}\xi_i(l) = \xi^{\textrm{IC}}_{i}(l),\quad i\in\mathcal{T}\label{P6}.
\end{align}
\end{subequations}

Constraint (\ref{P7}) are the pipe and actuation constraints described in Section~\ref{sec:act}. Constraint (\ref{P1}) is the convex relaxation of the growth kinetics; it is implemented in SOC form, e.g., as (\ref{eq:SOCC}) or (\ref{eq:SOCM}). Constraint (\ref{P2}) is the approximation of the discretized dynamics of the plant concentrations. Constraint (\ref{P3}) is the discretized dynamics of the volumes in the network. Constraints (\ref{P4})-(\ref{P6}) are the initial conditions of the trajectory.

The initial conditions in Constraints (\ref{P4})-(\ref{P6}) are the observations or estimates of the current state variables. E.g., If the current period is $n$, then $V^{\textrm{IC}}_{i}(-l)=\breve{V}_i(n-l)$ for $l=1,...,\tau^{\textrm{IC}}$. The other problem parameters, e.g., influent forecasts, should also reflect the actual system time. Once the optimization is solved, the control consists of the optimal actuator flow rates in the first time period:
\[
Q_{ij}(n)=Q^*_{ij}(0),\quad ij\in\mathcal{A}.
\]

\subsection{Conventional flow-rate controller}\label{sec:sewer:flowcontroller}

We now describe a conventional, volume-based MPC that minimizes flooding and CSO~\cite{ocampo2010model,ocampo2013applicationfull,garcia2015modeling}. It is the same as our controller except for the following changes to the optimization $\texttt{Traj}(\cdot)$.
\begin{itemize}
\item Omit all concentration and reaction variables, $\xi$ and $T$.
\item Omit the pollutant release and microbial growth terms from the objective.
\item Omit constraints (\ref{P0})-(\ref{P2}) and (\ref{P6}).
\end{itemize}
Because the volume-based controller's optimization does not contain concentrations or growth kinetics, it has only linear constraints and linear and convex quadratic terms in the objective. 

\section{Simulations}\label{sec:simulations}

\begin{figure}[h]
\centering
\includegraphics[width=\columnwidth]{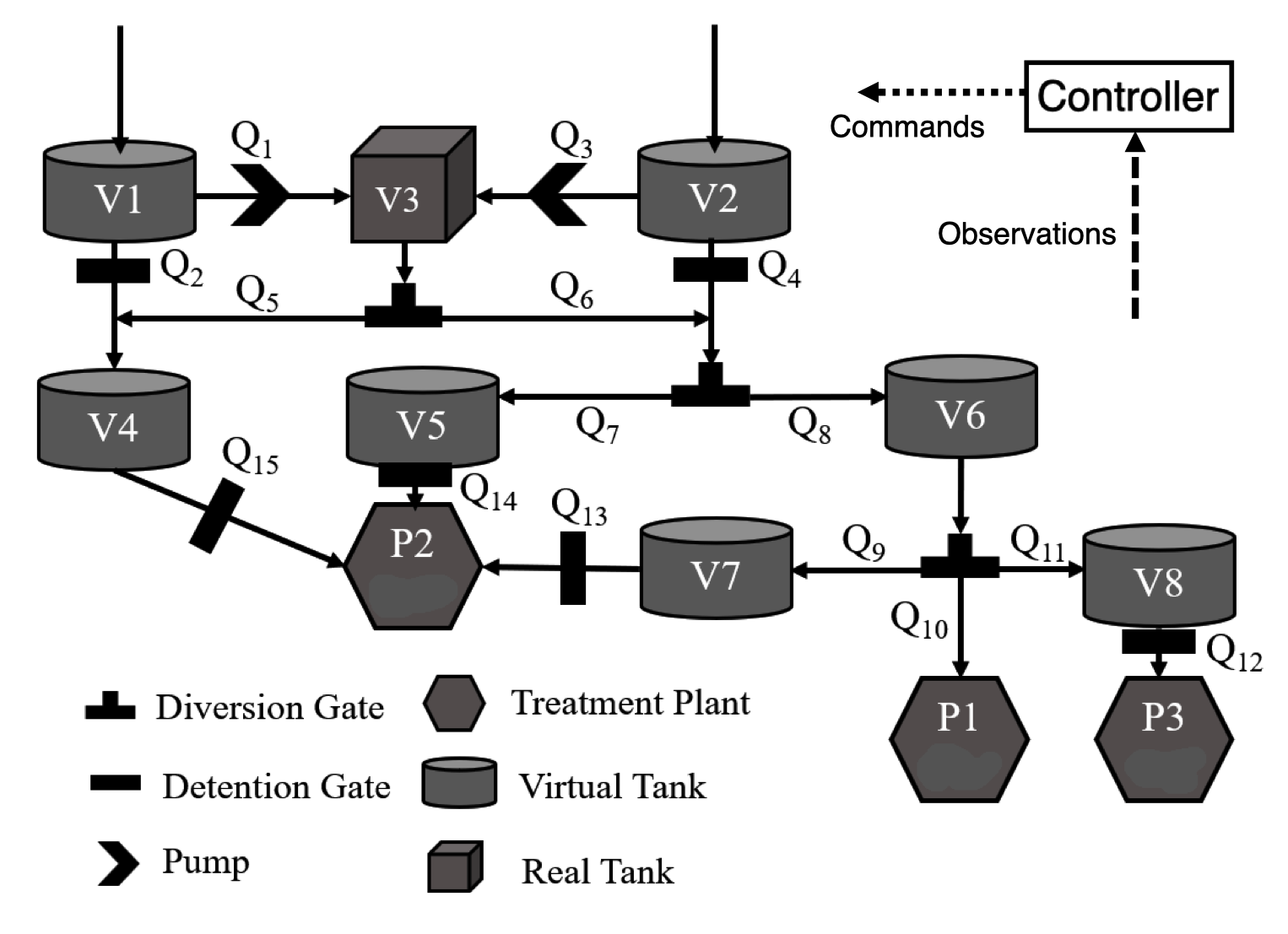}
\caption{Test system.}
\label{fig:paris}
\end{figure}

The test system, shown in Figure~\ref{fig:paris}, is based on a portion of the Paris wastewater treatment network. The network can store sewage in one storage tank and seven long sewer pipes, which are virtual tanks. There are two pumps, five detention gates, and three diversion gates, which regulate how sewage arrives at the three treatment plants. A complete system description is given in~\cite{taylor2024predictive}.

We model one biomass, $X$, and four substrates: biochemical oxygen demand (BOD), ammonia ($\textrm{NH}_4$), nitrite ($\textrm{NO}_2^-$), and nitrate ($\textrm{NO}_3^-$). There are four biochemical reactions, all with Contois growth rates. The parameters, shown Table~\ref{GrowthPar}, were taken from~\cite{ROBLESRODRIGUEZ2018880}. Because we are using the Contois and not Monod growth rate, we scaled the half-saturation constants by $10^{-3}$.

\begin{table}[h]
	\centering
	\begin{tabular}{|l|c|c|c|}
	\hline
	Parameter & Plant 1 & Plant 2 & Plant 3  \\
	\hline
	$\mu^{\textrm{BOD}}$ $\left(\textrm{day}^{-1}\right)$ & $3.99$ & $2.56$ & $1.93$\\
	$\mu^{\textrm{NH}_4^+}$ & $0.84$ & $0.83$ & $0.89$\\
	$\mu^{\textrm{NO}_2^-}$ & $1.68$ & $1.27$ & $0.92$\\
	$\mu^{\textrm{NO}_3^-}$ & $1.21$ & $1.38$ & $0.85$\\
	\hline
	$K^{\textrm{BOD}}\;(\times 10^3)$ & $13.67$ & $11.65$ & $14.26$\\
	$K^{\textrm{NH}_4^+}$ & $6.59$ & $14.98$ & $8.53$\\
	$K^{\textrm{NO}_2^-}$ & $2.46$ & $1.15$ & $2.55$\\
	$K^{\textrm{NO}_3^-}$ & $1.40$ & $2.69$ & $4.20$\\
	\hline
	$y^{\textrm{NH}_4^+,\textrm{NO}_2^-}$ & $0.28$ & $0.25$ & $0.27$\\
	$y^{\textrm{NO}_2^-,\textrm{NO}_3^-}$ & $0.68$ & $0.64$ & $0.70$\\
	$y^{X,\textrm{BOD}}$ & $0.67$ & $0.67$ & $0.67$\\
	$y^{X,\textrm{NH}_4^+}$ & $0.24$ & $0.24$ & $0.24$\\
	\hline
	\end{tabular}
	\caption{Reaction parameters}
	\label{GrowthPar}
\end{table}
The stoichiometric matrix for each plant $i\in\mathcal{S}$ is
\[
\kappa_i=\left[
\begin{array}{cccc}
-1&0 & 0 & 0\\
0 & -1 & 0 & 0\\
0 & 1/y_i^{\textrm{NH}_4^+,\textrm{NO}_2^-} & -1 & 0\\
0 & 0 & 1/y_i^{\textrm{NO}_2^-,\textrm{NO}_3^-} & -1\\
y_i^{X,\textrm{BOD}} & y_i^{X,\textrm{NH}_4^+} 0 &0  &0 \\
\end{array}
\right].
\]
The biomass death rate in Plants 1, 2, and 3 is 0.01, 0.1, and 0.1 $\textrm{day}^{-1}$, respectively, and the biomass portion of the effluent has a tenth the outflow rate of the substrates.

\begin{figure}[h]
\centering
\includegraphics[width=\columnwidth]{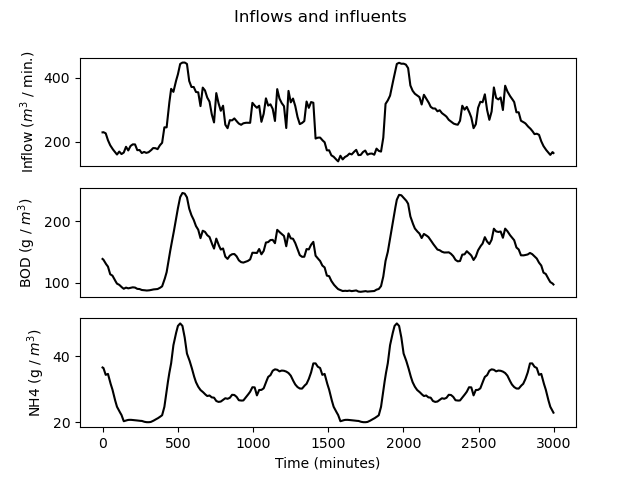}
\caption{Influents}
\label{fig:inflow}
\end{figure}

The inflow and influent data is from the \texttt{Inf\_dry\_2006} data set of~\cite{IWABSM1}, with the inflows scaled to be within system capacity. Volumes V1 and V2 receive identical influent. The flow rate and concentrations of BOD and $\textrm{NH}_4$ into each of the two volumes through time are shown in Figure~\ref{fig:inflow}. The influent concentrations of $\textrm{NO}_2^-$ and $\textrm{NO}_3^-$ are zero.

We compare our controller (\texttt{FC}) to the conventional volume-based controller (\texttt{F}) described in Section~\ref{sec:sewer:flowcontroller}. \texttt{FC} seeks to minimize pollutant release and maximize microbial growth. \texttt{F} seeks to minimize the stored volume, thus passing through the flow as quickly as possible. Both controllers promote smooth actuation via slope and curvature objectives. All optimizations were implemented using the parser CVXPy~\cite{diamond2016cvxpy} and solver Clarabel~\cite{goulart2024clarabel}. The plots were created with MatPlotLib~\cite{hunter2007matplotlib}.

The simulation is over 50 hours. Each control time period is 15 minutes long, so that there are 200 total time periods in the simulation. The trajectory optimizations of both controllers use the third-order Adams-Moulton discretization with a time step of $\Delta_{\textrm{T}}=3$ minutes long. The time horizon of both controllers is $H\Delta_{\textrm{T}}= 8$ hours, so that there are $H=160$ time periods in each trajectory optimization. The actuator variables are constrained to be constant within each control time period. \texttt{Traj}$(\cdot)$ is an SOCP with 21,091 variables and 42,807 constraints.

Table~\ref{tab:metrics} shows the key performance metrics for each controller. Solving the trajectory optimization in \texttt{FC} takes about three times as long as that in \texttt{F}. In both cases, the computation time is more than an order of magnitude less than the 15-minute duration of each control time period. Under both controllers, the system quickly hits its capacity limits due to the high initial volumes in the network, but avoids any flooding or CSO. Both controllers treat roughly the same volume of sewage over 50 hours. \texttt{FC} releases roughly 15.6\% less total pollutant mass than \texttt{F}.

\begin{table}[h]
	\centering
	\begin{tabular}{|l|l|l|}
	\hline
	 Metric& \texttt{F} & \texttt{FC}  \\
	\hline
	Mean trajectory optimization running time (s)& 15 & 44 \\
	Total treated volume ($10^6\times \textrm{m}^3$) & 9.30 & 9.23 \\
	Total pollutant release ($10^4\times$ kg) & 3.09 & 2.61 \\
	\hline
	\end{tabular}
 	\caption{Performance metrics}
	\label{tab:metrics}
\end{table}

\begin{figure}[h]
\centering
\includegraphics[width=\columnwidth]{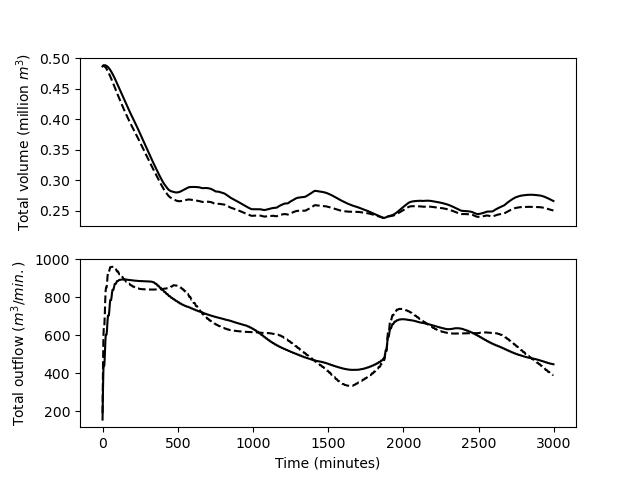}
\caption{Total volume and outflow}
\label{fig:total}
\end{figure}

\begin{figure}[h]
\centering
\includegraphics[width=\columnwidth]{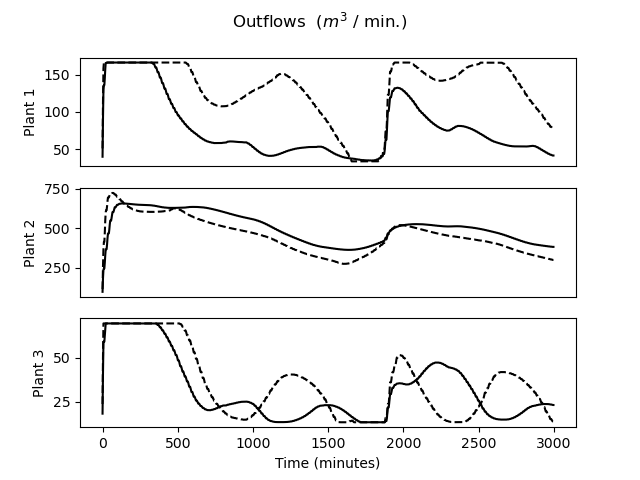}
\caption{Plant outflows}
\label{fig:outflow}
\end{figure}

\begin{figure}[h]
\centering
\includegraphics[width=\columnwidth]{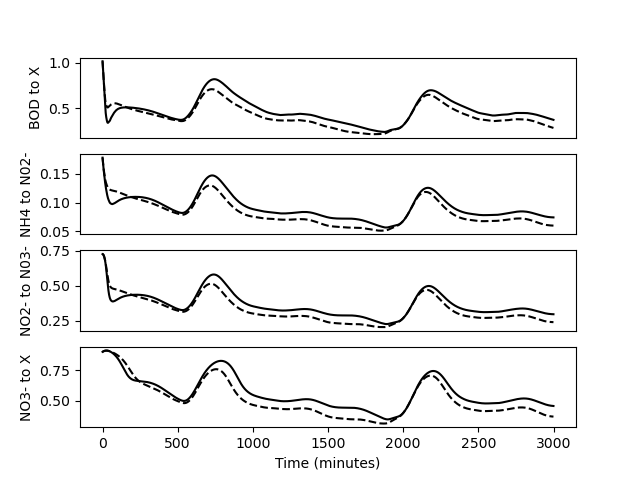}
\caption{Reaction rates in Plant 2}
\label{fig:reactionrates}
\end{figure}

\begin{figure}[h]
\centering
\includegraphics[width=\columnwidth]{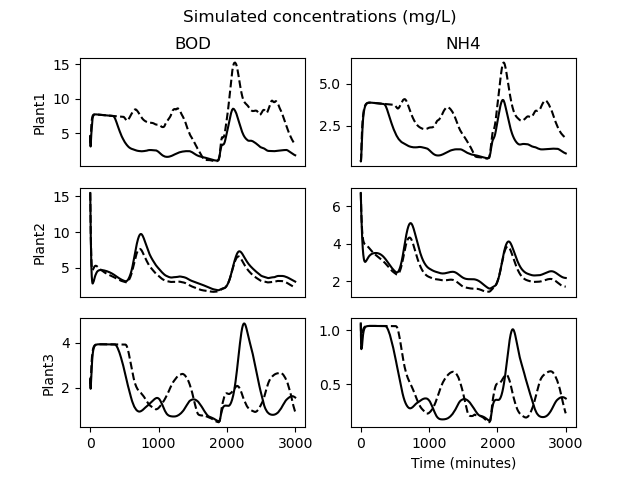}
\caption{Plant concentrations}
\label{fig:concentration}
\end{figure}

In Figures~\ref{fig:total} through \ref{fig:concentration}, the solid lines correspond to \texttt{FC} and the dashed lines to \texttt{F}. Figure~\ref{fig:total} shows that \texttt{FC} retains a higher total volume in the network than \texttt{F}, and discharges more evenly over time.
Figure~\ref{fig:outflow} shows the outflow from each of the three plants. \texttt{FC} passes more flow through Plant 2 and less through Plant 1 than the \texttt{F}. 
Figure~\ref{fig:reactionrates} shows that this leads to higher reaction rates in Plant 2, and thus more efficient treatment in the largest of the three plants.
\texttt{FC} avoids the concentration spike in Plant 1 around 2,000 minutes, as shown in Figure~\ref{fig:concentration} for BOD and $\textrm{NH}_4$. The same spike is also avoided for $\textrm{NO}_2^-$ and $\textrm{NO}_3^-$.

The broader trend in this case study is that \texttt{FC} better balances the flow over plants and over time, achieving better treatment through higher reaction rates and by avoiding large high-concentration discharges.

\section{Conclusion}

We have designed a novel MPC to minimize flooding and CSO during rain and to balance the flow across plants and time during dry weather, reducing pollutant release. It is the first pollution-based controller that relies entirely on convex optimization. The underlying optimization is nonconvex due to nonlinear microbial growth and bilinear mass flows. We have convexified the problem via relaxation of the growth kinetics and physically motivated linearization of the bilinearities. The resulting trajectory optimization in each time period is an SOCP, which an be quickly solved at large scales. In simulation, the new controller releases 15\% less pollutant mass while treating nearly the same volume as a conventional volume-based MPC.

\bibliographystyle{IEEEtran}
\bibliography{/Users/joshtaylor/Dropbox/Work/Projects/MainBib,/Users/joshtaylor/Dropbox/Work/Projects/JATBib}
\end{document}